\pdfoutput=1
\documentclass[a4paper]{amsart}
\usepackage[T1]{fontenc}
\usepackage{amsthm}
\usepackage{amssymb}
\usepackage[unicode=true,pdfusetitle,
 bookmarks=true,bookmarksnumbered=false,bookmarksopen=false,
 breaklinks=false,pdfborder={0 0 0},pdfborderstyle={},backref=false,colorlinks=false]
 {hyperref}
\theoremstyle{plain}
\newtheorem{thm}{Theorem}[section]
\theoremstyle{definition}

\theoremstyle{remark}

\theoremstyle{remark}

\theoremstyle{definition}

\theoremstyle{definition}
\newtheorem{construction}[thm]{Construction}
\theoremstyle{plain}
\newtheorem{prop}[thm]{Proposition}
\theoremstyle{plain}

\theoremstyle{plain}

\theoremstyle{plain}

\usepackage{tikz-cd}
\usepackage{mathtools}
\usepackage{adjustbox}
\usepackage{enumitem}
\tikzcdset{scale cd/.style={every label/.append style={scale=#1},
    cells={nodes={scale=#1}}}}
\usepackage{eucal}
\usepackage{mleftright}
\mleftright

\newcommand{\notehelper}[3]{\textcolor{#3}{$\blacksquare$}\marginpar{\ifodd\thepage\raggedright\else\raggedleft\fi\color{#3}\tiny \textbf{#2:} #1}}

\usetikzlibrary{calc}
\usetikzlibrary{decorations.pathmorphing}
\tikzset{curve/.style={settings={#1},to path={(\tikztostart)
			.. controls ($(\tikztostart)!\pv{pos}!(\tikztotarget)!\pv{height}!270:(\tikztotarget)$)
			and ($(\tikztostart)!1-\pv{pos}!(\tikztotarget)!\pv{height}!270:(\tikztotarget)$)
			.. (\tikztotarget)\tikztonodes}},
	settings/.code={\tikzset{quiver/.cd,#1}
		\def\pv##1{\pgfkeysvalueof{/tikz/quiver/##1}}},
	quiver/.cd,pos/.initial=0.35,height/.initial=0}

\global\long\def\sf#1{\mathsf{#1}}%
\global\long\def\cal#1{\mathcal{#1}}%
\global\long\def\pr#1{\left(#1\right)}%
\global\long\def\abs#1{\left|#1\right|}%
\global\long\def\Catinfty{\sf{Cat}_{\infty}}%
\global\long\def\An{\sf{An}}%
\global\long\def\RC{\sf{RelCat}_{\infty}}%
\global\long\def\lim{\operatorname{lim}}%
\global\long\def\colim{\operatorname{colim}}%
\global\long\def\Fun{\operatorname{Fun}}%
\global\long\def\Map{\operatorname{Map}}%
\global\long\def\op{\mathrm{op}}%
\global\long\def\id{\mathrm{id}}%
\global\long\def\ac{\mathrm{ac}}%
\global\long\def\Del{\mathbf{\Delta}}%
\global\long\def\rel{\mathrm{rel}}%
\newcommand{\iso}{\xrightarrow{\;\smash{\raisebox{-0.5ex}{\ensuremath{\scriptstyle\sim}}}\;}}

\title[A short proof of the universality of the relative Rezk nerve]{A short proof of the universality \\ of the relative Rezk nerve}
\author{Kensuke Arakawa}
\address{K.A.: Department of Mathematics, Kyoto University, 6068502, Kyoto, Japan}
\author{Bastiaan Cnossen}
\address{B.C.: Fakultät für Mathematik, Universität Regensburg, 93040 Regensburg, Germany}

\begin{document}
\begin{abstract}
    We give a concise, conceptual proof of the universality of the relative Rezk
    nerve, due to Mazel-Gee.
\end{abstract}
\subjclass[2020]{18N55, 18N60, 18A22}
\maketitle
\section{Introduction}
Localization, the process of formally inverting a class of morphisms in a category, is a foundational concept in category theory. 
A primary example is the derived category of an abelian category, obtained by inverting the quasi-isomorphisms in the category of chain complexes.

The definition of localization can easily be generalized to the $\infty$-categorical setting. Starting from an $\infty$-category $\mathcal{C}$ and a collection $\mathcal{W}$ of morphisms in $\mathcal{C}$, one can construct a new $\infty$-category $\mathcal{C}[\mathcal{W}^{-1}]$ characterized up to equivalence by the property that functors out of it correspond to functors out of $\mathcal{C}$ which map morphisms in $\mathcal{W}$ to equivalences. It turns out that this generalization is at least as important as localizations of ordinary categories, because
many $\infty$-categories of interest arise as $\infty$-categorical localizations. For example, the full $\infty$-category $\An \subseteq \Catinfty$ spanned by the \textit{animae}/\textit{$\infty$-groupoids} can be obtained from the category of topological spaces by inverting the weak homotopy equivalences. However, the universal property defining $\mathcal{C}[\mathcal{W}^{-1}]$ is abstract and provides little direct insight into the structure of the localized $\infty$-category itself, such as its mapping spaces. A fundamental challenge, therefore, is to find more explicit constructions for $\mathcal{C}[\mathcal{W}^{-1}]$. This challenge has driven intensive research since the early days of $\infty$-category theory, see for instance \cite{DK80_1, DK80_2, DK80_3}, \cite[$\S$5.2.7 and $\S$6]{HTT}, and \cite[Chapter 7]{Cisinski19}.

In this paper, we will revisit an elegant presentation of $\mathcal{C}[\mathcal{W}^{-1}]$ due to Mazel-Gee \cite{MG19}. To state it, recall that a powerful strategy for understanding and constructing $\infty$-categories involves ``nerve'' constructions, which translate $\infty$-categorical data into the world of simplicial objects. A primary example is the \textbf{Rezk nerve functor}
\[
N\colon\Catinfty \hookrightarrow \Fun\pr{\Del^{\op},\An}, \qquad N(\cal C)_n := \Map_{\Catinfty}\pr{[n],\cal C},
\]
which is a fully faithful embedding of the $\infty$-category of (small) $\infty$-categories into that of simplicial animae; see \cite{Rezk01, HS25}. This construction can be extended to \textbf{relative $\infty$-categories}, i.e., pairs $\pr{\cal C,\cal W}$, where $\cal C$ is an $\infty$-category and $\cal W\subset\cal C$ is a subcategory containing all equivalences, whose morphisms are called \textbf{weak equivalences}. By regarding relative $\infty$-categories equivalently as the essentially surjective monomorphisms of $\infty$-categories $\cal W \hookrightarrow \cal C$, they define a full subcategory $\RC \subseteq \Fun\pr{[1],\Catinfty}$, and we define the \textbf{relative Rezk nerve functor} as
\[
N^{\rel}\colon\RC \to\Fun\pr{\Del^{\op},\An}, \qquad N^{\rel}(\cal C, \cal W)_n := \abs{\Fun\pr{[n],\cal C}\times_{\cal C^{n+1}}\cal W^{n+1}},
\]
where the realization functor $\abs -\colon\Catinfty\to\An$ is the left adjoint to the inclusion $\An \hookrightarrow \Catinfty$. 

In case $\mathcal{W}$ is the core of $\mathcal{C}$, the relative nerve recovers the ordinary nerve, i.e., there is a natural equivalence $N^{\rel}(\mathcal{C}, \mathcal{C}^{\simeq}) \simeq N(\mathcal{C})$. The inclusion of relative $\infty$-categories $(\mathcal{C},\mathcal{C}^{\simeq}) \hookrightarrow (\mathcal{C},\mathcal{W})$ thus induces a natural comparison map
\[
\mathcal{C} \simeq \ac(N(\mathcal{C})) \to \ac(N^{\rel}(\mathcal{C},\mathcal{W})),
\]
where we write $\ac\colon\Fun\pr{\Del^{\op},\An}\to\Catinfty$ for the left adjoint to the (fully faithful) Rezk nerve functor. Mazel-Gee's theorem states that this map presents the associated $\infty$-category of the relative nerve as a model for the localization:

\begin{thm}[{\cite[Theorem~3.8]{MG19}}]\label{thm:MG}
	The map $\mathcal{C} \to \ac\pr{N^{\rel}\pr{\cal C,\cal W}}$ exhibits its target as the localization of $\mathcal{C}$ at $\mathcal{W}$, i.e., it induces an equivalence
	\[
	\cal C[\cal W^{-1}] \iso \ac\pr{N^{\rel}\pr{\cal C,\cal W}}.
	\]
\end{thm}

The relative Rezk nerve of relative 1-categories was introduced by Rezk in \cite[Section~3]{Rezk01}, where he called it the \textbf{classification diagram functor}. Its relation to localization had been partially known prior to Mazel-Gee's work. Indeed, already when Rezk introduced the classification diagram, he proved a special case of Theorem \ref{thm:MG} where $\mathcal{C}$ is assumed to be the underlying category of a simplicial model category and $\mathcal{W}$ the subcategory of weak equivalences. Bergner then generalized this to arbitrary model categories \cite[Theorem 6.2]{Ber09}, further deepening the understanding of the relation. The classification diagram functor also played a key role in Barwick--Kan's papers \cite{BK12, BK12b}, where they showed that relative 1-categories model $\infty$-categories via localization. Mazel-Gee completed this picture by establishing Theorem \ref{thm:MG} in the full general setting of $\infty$-categories.

Mazel-Gee's theorem is a powerful tool for accessing and understanding localizations, and it has found numerous applications;
see \cite{LNS17, HHFSN23, KSW24, A24b, Hinich24} for instance.
However, Mazel-Gee's original proof of the theorem relies on a rather intricate point-set argument, which he himself acknowledges as ‘unsatisfying.’
While the first author later simplified the proof (and proved a generalization) in \cite{A23b}, this simplification still depended on specific models of $\infty$-categories.
In this paper, we address this situation by presenting a short proof of Theorem \ref{thm:MG} that relies only on basic, intrinsic constructions within the theory of $\infty$-categories.

\subsection*{Notation and conventions}
Given relative $\infty$-categories $\pr{\cal C,\cal W}$ and $\pr{\cal C',\cal W'}$, we write
\[
\Fun_{\rel}\pr{\pr{\cal C,\cal W},\pr{\cal C',\cal W'}}\subset\Fun\pr{\cal C,\cal C'}
\]
for the full subcategory spanned by the functors $\cal C\to\cal C'$ carrying morphisms in $\cal W$ to morphisms in $\cal W'$. We identify $\Catinfty$ with a full subcategory of $\RC$ via the embedding $\cal C\mapsto\pr{\cal C,\cal C^{\simeq}}$, where $\cal C^{\simeq}$ denotes the core of $\cal C$.

	\section{Proof of Theorem \ref{thm:MG}}
	
	The main ingredient for our proof of Theorem \ref{thm:MG} is a simplicial resolution of localizations.
	\begin{construction}
		\label{const:sRes}We define a functor $L_{\bullet}\colon\RC\to\Fun\pr{\Del^{\op},\Catinfty}$
		by the formula
		\[
		L_{n}\pr{\cal C,\cal W} := \Fun_{\rel}\pr{\pr{[n],[n]},\pr{\cal C,\cal W}}.
		\]
		Thus $L_{0}\pr{\cal C,\cal W}$ identifies with $\cal C$,
		$L_{1}\pr{\cal C,\cal W}\subset\Fun\pr{[1],\cal C}$ is the full subcategory
		of weak equivalences (not to be confused with $\Fun\pr{[1],\cal W}$), and so on. 
	\end{construction}
	
	Since for all $\cal D \in \Catinfty$ and $[n] \in \Del^{\op}$ the constant diagram functor $\cal D \to L_n(\cal D) := L_n(\cal D, \cal D^{\simeq})$ is an equivalence, the localization functor $\gamma\colon \cal C \to \cal C [\cal W^{-1}]$ induces maps
	\begin{equation}\label{eq:sRes}
		L_n(\gamma)\colon L_n(\cal C, \cal W) \to L_n(\cal C[\cal W^{-1}]) \simeq \cal C[\cal W^{-1}]
	\end{equation}
	naturally in $[n] \in \Del^{\op}$ and $(\cal C, \cal W) \in \RC$.
		
	\begin{prop}\label{lem:sRes}
		For a relative $\infty$-category $\pr{\cal C,\cal W}\in\RC$, the functor
		\[
			\colim_{\Del^{\op}}L_{\bullet}\pr{\cal C,\cal W} \to \cal C[\cal W^{-1}]
		\]
		induced by \eqref{eq:sRes} is an equivalence of $\infty$-categories, natural in $(\cal C, \cal W) \in \RC$.
	\end{prop}
	\begin{proof}
		The functor in question is the composite
		\[
			\colim_{\Del^{\op}}L_{\bullet}\pr{\cal C,\cal W} \xrightarrow{\colim_{\Del^{\op}} L_{\bullet}(\gamma)} \colim_{\Del^{\op}}L_{\bullet}\pr{\cal C[\cal W^{-1}]} \iso \cal C[\cal W^{-1}].
		\]
		We may thus equivalently show that for each $\infty$-category $\cal D$, the induced map
		\[
		\theta\colon\Map_{\Fun\pr{\Del^{\op},\Catinfty}}\pr{L_{\bullet}\pr{\cal C[\cal W^{-1}]},\delta\pr{\cal D}}\to\Map_{\Fun\pr{\Del^{\op},\Catinfty}}\pr{L_{\bullet}\pr{\cal C,\cal W},\delta\pr{\cal D}}
		\]
		is an equivalence, where $\delta\pr{\cal D}$ denotes the constant simplicial object at $\cal D$. 
		
		\textit{Step 1:} We first show that $\theta$ is an inclusion of a full subanima, whose image consists of those maps $\phi\colon L_{\bullet}\pr{\cal C,\cal W}\to\delta\pr{\cal D}$
		with the following property:
		\begin{itemize}
			\item [($\ast$)]For each $n\geq0$, the map $\phi_{n}\colon\Fun_{\rel}\pr{\pr{[n],[n]},\pr{\cal C,\cal W}}\to\cal D$
			carries natural weak equivalences to equivalences in $\cal D$.
		\end{itemize}
		By applying \cite[Proposition A.11]{Ram23} to $\Catinfty^{\op}$ and using that localizations of $\infty$-categories are monomorphisms in $\Catinfty^{\op}$, this claim reduces to showing that for each $n$ the map $L_{n}(\gamma)\colon L_{n}\pr{\cal C,\cal W} \to L_{n}\pr{\cal C[\cal W^{-1}]}$ is a localization at the natural weak equivalences. To this end, consider the commutative diagram
		\[
		\begin{tikzcd}
			\cal C \simeq L_{0}\pr{\cal C,\cal W} \rar{L_0(\gamma)} \dar & L_0(\cal C[\cal W^{-1}]) \simeq \cal C[\cal W^{-1}] \dar[shift right= 8]{\sim} \\
			L_{n}\pr{\cal C,\cal W} \rar{L_n(\gamma)} & L_n(\cal C[\cal W^{-1}]) \simeq \cal C[\cal W^{-1}]
		\end{tikzcd}
		\]
		obtained by naturality for $[n] \to [0]$. Since the right vertical map is an equivalence and $L_0(\gamma) = \gamma$ is a localization at $\cal W$, it remains to show that the left vertical map $\cal C \to \operatorname{Fun}(\pr{[n],[n]},(\mathcal{C},\mathcal{W}))$ induces an equivalence upon localizing at the natural weak equivalences. But this is clear, because the evaluation at $0\in[n]$ is a right adjoint and the unit and counit transformations are natural weak equivalences. 
		
		\textit{Step 2:} For surjectivity of $\theta$, we must show that every map $\phi\colon L_{\bullet}\pr{\cal C,\cal W}\to\delta\pr{\cal D}$
		satisfies condition ($\ast$). Using the naturality of $\phi_{n}$ in $[n]$, it suffices to prove ($\ast$) for $n=0$. Let $w\colon X\to Y$
		be a weak equivalence of $\cal C$. We wish to show that $\phi_{0}\pr w$
		is an equivalence. There is a morphism $\alpha\colon w\to\id_{Y}$ in $L_1(\cal C, \cal W) = \Fun_{\rel}\pr{\pr{[1],[1]},\pr{\cal C,\cal W}}$
		depicted as
		\[\begin{tikzcd}
			X & Y \\
			Y & Y.
			\arrow["w", from=1-1, to=1-2]
			\arrow["w"', from=1-1, to=2-1]
			\arrow["{\operatorname{id}_Y}", from=1-2, to=2-2]
			\arrow["{\operatorname{id}_Y}"', from=2-1, to=2-2]
		\end{tikzcd}\]
		Using the naturality of the maps $\{\phi_{n}\}_{[n]\in\Del^{\op}}$,
		we obtain
		\[
		\phi_{0}\pr w\simeq\phi_{0}\pr{d_{1}\pr{\alpha}}\simeq\phi_{1}\pr{\alpha}\simeq\phi_{0}\pr{d_{0}\pr{\alpha}}\simeq\phi_{0}\pr{\id_{Y}}.
		\]
		The right-hand side is evidently an equivalence, so $\phi_{0}\pr w$
		is also an equivalence, as desired.
	\end{proof}
	\begin{proof}
		[Proof of Theorem \ref{thm:MG}] We contemplate the following commutative diagram:
		\[\hspace{-7pt}\begin{tikzcd}[column sep=tiny]
			{\mathsf{RelCat}_\infty} && {\operatorname{Fun}(\mathbf{\Delta}^{\mathrm{op}},\mathsf{Cat}_\infty)} \\
			{\operatorname{Fun}(\mathbf{\Delta}^{\mathrm{op}},\mathsf{Cat}_\infty)} \\
			& {\operatorname{Fun}(\mathbf{\Delta}_1^{\mathrm{op}},\operatorname{Fun}(\mathbf{\Delta}_2^{\mathrm{op}},\mathsf{An}))} \\
			{\operatorname{Fun}(\mathbf{\Delta}^{\mathrm{op}},\mathsf{An})} && {\operatorname{Fun}(\mathbf{\Delta}_2^{\mathrm{op}},\operatorname{Fun}(\mathbf{\Delta}_1^{\mathrm{op}},\mathsf{An}))} \\
			{\mathsf{Cat}_\infty} && {\operatorname{Fun}(\mathbf{\Delta}^{\mathrm{op}},\mathsf{Cat}_\infty)}
			\arrow["{L_{\bullet}}", from=1-1, to=1-3]
			\arrow["F", from=1-1, to=2-1]
			\arrow[""{name=0, anchor=center, inner sep=0}, "{N^{\mathrm{rel}}}"', shift right=5, curve={height=30pt}, from=1-1, to=4-1]
			\arrow["{(3)}"{description, pos=0.6}, draw=none, from=1-3, to=3-2]
			\arrow["{\operatorname{Fun}(\mathbf{\Delta}^{\mathrm{op}},N)}"{description}, from=1-3, to=4-3]
			\arrow[""{name=1, anchor=center, inner sep=0}, "{\operatorname{id}}", shift left=15, curve={height=-30pt}, from=1-3, to=5-3]
			\arrow["{\operatorname{Fun}(\mathbf{\Delta}^{\mathrm{op}},N)}"{description}, from=2-1, to=3-2]
			\arrow[""{name=2, anchor=center, inner sep=0}, "{\operatorname{Fun}(\mathbf{\Delta}^{\mathrm{op}},|-|)}"{description}, from=2-1, to=4-1]
			\arrow["{\operatorname{Fun}(\mathbf{\Delta}^{\mathrm{op}},\operatorname{colim})}"{description}, from=3-2, to=4-1]
			\arrow["\sigma", from=3-2, to=4-3]
			\arrow["{\mathrm{ac}}"', from=4-1, to=5-1]
			\arrow[""{name=3, anchor=center, inner sep=0}, "{\operatorname{colim}}"{description}, from=4-3, to=4-1]
			\arrow["{\operatorname{Fun}(\mathbf{\Delta}^{\mathrm{op}},\mathrm{ac})}"{description}, from=4-3, to=5-3]
			\arrow[""{name=4, anchor=center, inner sep=0}, "{\operatorname{colim}}", from=5-3, to=5-1]
			\arrow["{(4)}"{description, pos=0.3}, draw=none, from=0, to=2]
			\arrow["{(1)}"{description}, draw=none, from=3-2, to=3]
			\arrow["{(6)}"{description, pos=0.2}, draw=none, from=3-2, to=2]
			\arrow["{(2)}"{description}, draw=none, from=3, to=4]
			\arrow["{(5)}", shift left=3, draw=none, from=4-3, to=1]
		\end{tikzcd}\]
		Here $F$ is given by $F\pr{\cal C,\cal W}_{n}=\Fun\pr{[n],\cal C}\times_{\cal C^{n+1}}\cal W^{n+1}$,
		$\Del_{1}$ and $\Del_{2}$ are copies of $\Del$, and $\sigma$ interchanges
		the bisimplicial degrees. The commutativity of (1) and (2) is clear, (3) and (4) commute by unwinding definitions, and (5) commutes by fully faithfulness of $N$. Commutativity of (6) follows from the natural equivalence $\abs{\cal C} \simeq \abs{\ac N \cal C} \simeq \colim N \cal C$, where the first equivalence again uses fully faithfulness of $N$ and the equivalence $\abs{-} \circ \ac \simeq \colim$ is clear upon passing to right adjoints. 
		
		The outer square in the diagram thus provides a natural equivalence $\ac N^{\rel}(\cal C, \cal W) \simeq \colim_{\Del^{\op}} L_\bullet(\cal C, \cal W)$. Theorem \ref{thm:MG} consequently reduces to the claim that the map $\cal C \simeq \colim_{\Del^{\op}} L_\bullet(\cal C) \to \colim_{\Del^{\op}} L_\bullet(\cal C, \cal W)$ exhibits the target as the localization of $\cal C$ at $\cal W$, which is the content of Proposition \ref{lem:sRes}.
	\end{proof}

\subsection*{Acknowledgment}
K.A. was supported by JSPS KAKENHI Grant Number 24KJ1443. B.C. is an associate member of the SFB 1085 Higher Invariants.
\providecommand{\bysame}{\leavevmode\hbox to3em{\hrulefill}\thinspace}
\providecommand{\MR}{\relax\ifhmode\unskip\space\fi MR }
\providecommand{\MRhref}[2]{%
  \href{http://www.ams.org/mathscinet-getitem?mr=#1}{#2}
}
\providecommand{\href}[2]{#2}

\end{document}